\newcount\hh
\newcount\mm
\mm=\time
\hh=\time
\divide\hh by 60
\divide\mm by 60
\multiply\mm by 60
\mm=-\mm
\advance\mm by \time
\def\hhmm{\number\hh:\ifnum\mm<10{}0\fi\number\mm}
\documentclass[twoside]{article}
\usepackage{amssymb, amsthm,amsmath}
\usepackage{pictexwd, dcpic}
\usepackage[]{hyperref}
\hypersetup{
colorlinks=true,
linkcolor=black,
}
\usepackage{a4wide} % improves hbox error and tolerance on width
\usepackage{fancyhdr}
\theoremstyle{definition} 
\newtheorem{theorem}{Theorem}
\newtheorem{prop}[theorem]{Proposition}
\newtheorem{lemma}[theorem]{Lemma}
\newtheorem*{defi}{Definition}
\newtheorem*{ack}{Acknowledgement}
\newtheorem*{exam}{Example}
\newtheorem{cor}[theorem]{Corollary}
\newtheorem{remark}[theorem]{Remark}

\theoremstyle{remark} 
\newtheorem*{rem}{Remark} %%% * means no numbering
\newtheorem*{nota}{Notation}
\newcommand\ann{\mathrm{Ann}}
\newcommand{\comment}[1]{}

\newcommand\brac[1]{\langle #1 \rangle}
\newcommand\sep{\,:\,}

\newcommand\st{\,\backepsilon\,}
\newcommand\qf[1]{\mathrm{Quot}(#1)}	
\newcommand\spec{\mathrm{Spec}\,}	
\newcommand\sper{\mathrm{Sper}\,}	
\newcommand\supp{\mathrm{supp}}	
\renewcommand\O{\mathcal{O}}

\newcommand\Z{\mathbb{Z}}

\newcommand\R{\mathbb{R}}

\newcommand{\N}{\mathbb{N}}
\newcommand{\p}{\mathfrak{p}}
\newcommand\wo[1]{\backslash{\{#1\}}}

\begin{document}

%%%%%%%%%%%%%%% TITLE %%%%%%%%
\title{Real Closed Rings and Real Closed * Rings}
\makeatletter
\let\mytitle\@title
\makeatother
\author{Jose Capco\\
\href{mailto:capco@fim.uni-passau.de}{\small{capco@fim.uni-passau.de}}\\
\emph{\small{Universit\"at Passau, Innstr. 33, 94032 Passau, Germany}} }
%\address{Universit\"at Passau, Innstr. 33, 94032 Passau, Germany}
%, Niels Schwartz}
%\email{capco@fim.uni-passau.de}
%\thanks{The author would like to thank the valuable advises of Prof. Schwartz and Marcus Tressl.}
%\subjclass{13J25}
\date{}
\thispagestyle{empty}
\maketitle
%%%%%%%%%%%%%%%%%%%%%%%%%%%

%%%%%%%%%%% Page Format %%%%%%%%%%%%%%%%%%%%%%%%%%
%\pagestyle{myheadings}
\pagestyle{fancy}
\fancyhead[R]{\mytitle}
\fancyhead[L]{J. Capco}
\tolerance=500
%%%%%%%%%%%%%%%%%%%%%%%%%%%%%%%%%%%%%%%%%%%%%%%%%%%

\begin{abstract}
Here we try to distinguish and compare different notions of real closedness mainly one developed by N. Schwartz in his 
Habilitationschrift and the other developed by A. Sankaran and K. Varadarajan in \cite{SV} which we shall call real closed *.  
We stick to the definition of real closed rings as defined and characterized in \cite{rcr} and we try to determine and 
characterize real closed rings that are real closed *. The main result is that real closed rings have unique real closure * 
and that real closure of real closed * rings arent necessarily real closed *. 
\comment{
We are concerned with the uniqueness of real 
closure * and we show that $f$-rings possesses this property.
}
\begin{description}
\item[Mathematics Subject Classification (2000):] Primary 13J25; Secondary 13B22, 13A35
\item[Keywords:] real closed rings, real closed $*$ rings, tight/essential extension, algebraic extension, total integral closure, 
Baer Hull, $f$-ring.
\end{description}
\end{abstract}

%\footnote{2000 Mathematics Subject Classification: Primary 13J25 ; Secondary 13B22, 13A35. Supported by DFG.}
\footnote{Supported by Deutsche Forschungsgemeinschaft and Universit\"at Passau.}
Unless otherwise stated all the rings in this paper are commutative unitary partially ordered rings (porings). The reader is assumed to know 
results in \cite{KS} Kapitel III.

\begin{nota}
Given a ring $A$ we use $\sper A$ to mean the topological space (Harrison Topology) of all the prime cones of $A$ such that they 
contain the partial ordering of $A$. If the partial ordering is not specifically mentioned then we will mostly assume it to be 
$$\sum A^2:=\{\sum_{i=1}^n a_i^2 \sep a_1,\dots,a_n\in A\}$$
For $a_1,\dots,a_n,a\in A$ we use the following notations:
$$\supp_A :\sper A \longrightarrow \spec A\qquad \supp_A(\alpha)=-\alpha\cap\alpha\quad \alpha\in\sper A $$
and 
\begin{eqnarray*}
V_A(a) & := & \{\p\in\spec A \sep a\in \p \} \\
Z_A(a) & := & \{\alpha\in\sper A \sep a\in\supp_A(\alpha)\} \\
%Z_A(a_1,\dots,a_n) & := & \{\alpha\in\sper A \sep a_1,\dots,a_n\in\supp_A(\alpha)\} \\
\bar P_A(a_1,\dots,a_n) & := & \{\alpha\in\sper A \sep a_1,\dots,a_n\in \alpha\} \\
\end{eqnarray*}
If it is clear with which ring we are working with, we drop the subscript $A$ above and we just write $\supp,V(a),\dots$ etc.

\end{nota}

\begin{defi}
\item[(a)] Let $A$ be a partially ordered integral domain, with partial ordering $A^+$, then we say that $A$ is a \emph{real closed
integral domain} iff the following holds
\begin{enumerate}
\item $\qf A$ is a real closed field
\item $A$ is integrally closed in $\qf A$
\item For any $a\in A$, the set $i(\sqrt{aA})$ is convex in $A_a$, where 
$i:A\rightarrow A_a$ is the canonical homomorphism.
\end{enumerate}
\item[(b)] A ring $A$ is a \emph{real closed ring} iff the following holds
\begin{enumerate}
\item $A$ is reduced
\item $\supp:\sper A\rightarrow \spec A$ is a homeomorphism
\item For any $\p\in\spec A$, $A/\p$ is real closed integral domain
\item If $a,b\in A$, then $\sqrt{aA} + \sqrt{bA}=\sqrt{aA+bA}$
\end{enumerate}
\end{defi}

\begin{remark}\label{rcrprop}
Most of the properties of real closed rings are discussed in \cite{Schw1} Chapter 1. Real closed rings have a lot of desirable algebraic
properties, the category has for instance limits and colimits (for small diagrams) thus fiber products of real closed rings are real 
closed as well. Taking factor rings, quotient rings and complete ring of quotients of a real closed ring will make it real 
closed as well (all have been worked out by N. Schwartz in \cite{Schw1,Gab,rcr,SM}).
\end{remark}

\begin{rem}
One sees that for a field this correspond to the usual definition of a real closed field.
\end{rem}

\begin{defi}
\begin{itemize}
\item Let $A$ be a ring, then a ring $B\supset A$ is said to be an \emph{essential extension} of $A$ iff for any ideal $I$ of $B$,
such that $I\neq {0}$ we have $I\cap A\neq {0}$. This is equivalent to saying that if $b\in B\wo 0$ then there is a $c\in B$ 
with $bc\in A\wo 0$.
\item A partially ordered ring $A$ with partial ordering $A^+$ is said to have bounded inversion iff $1+A^+ \subset A^*$
where $A^*$ is the set of units of $A$.
\end{itemize}
\end{defi}

\begin{defi}
A ring $A$ is real closed {*} (i.e. real closed according to \cite{SV}) iff there is no formally real ring that is an integral
extension and an essential extension of $A$.
\end{defi}

\begin{nota}
Give a reduced poring $A$, we use $\rho(A)$ to mean \underline{the} \emph{real closure of $A$}. The construction of the real closure of a reduced
poring can be read in \cite{SM} \S12. 

Also if $\alpha\in\sper A$, by $\rho(\alpha)$ we mean the real closure of the qoutient fields of 
$A/\supp(\alpha)$. 
\end{nota}

\begin{rem}
In \cite{Schw1} Chapter 2 and \cite{SM} \S12 the construction of the real closure of a partially ordered
ring is carried out. This construction is an extension of the ring itself and has the property that it is real closed and
that its real spectrum (and thus its prime spectrum as well) 
is canonically isomorphic to the real spectrum of the original ring.
\end{rem}

\begin{exam} We will show that there is a real closed {*} ring that isn't a real
closed ring. According to \cite{SV} Proposition 2, for domains, real
closed {*} is equivalent to the domain being integrally closed in
its qoutient field and that its qoutient field being a real closed
field. Now consider the ring $\Z$ and a real closed non-archimedian field $R$. According
to the Chevalley's Extension Theorem (Theorem 3.1.1 \cite{EP}) there is a valuation ring $V$ of $R$ with maximal ideal
$m_{V}$ such that $m_{V}\cap\Z=(2)$, and thus $V$ is not convex 
(because $2\in m_{V}$). However, we observe that $V$ is obviously
integrally closed in $R$ and that its qoutient field is $R$, so
$V$ is indeed real closed {*} and yet not real closed, since one sees from the the definition of real closed integral domains
that for $V$ to be real closed it is necessary that it be convex in $R$ and yet $2\in m_V$ which is a contradiction (see for instance
\cite{KS} Kapitel III, §2 Satz 3).
\end{exam}

\begin{nota}
For any commutative ring $A$ and any set $S\subset A$, we denote the annihilator of $S$ by 
$$\ann(S):=\{a\in A\sep aS=\{0\}\}$$
if $S$ is an $n$-tuple i.e. if $S=\{a_1,\dots,a_n\}\subset A$, then we may sometimes write $\ann(a_1,\dots,a_n)$ instead
of $\ann(\{a_1,\dots,a_n\})$ for simplicity.
\end{nota}

% BAER RINGS
\begin{defi}
Let $R$ be a commutative ring, then $R$ is said to be \emph{Baer} iff for any $S\subset R$ we can write
$$\ann(S)=eR$$
for some $e\in R$ such that $e^2=e$.
\end{defi}

We have the following Proposition that is obtained from \cite{Mewborn}

\begin{prop}\label{mborn}(Mewborn, 1971)
Let $R$ be a commutative reduced ring, then there is a ring $B(R)$ between $R$ and its complete ring of quotients $Q(R)$
such that $B(R)$ is Baer and $B(R)$ is contained in all Baer rings between $R$ and $Q(R)$ and in fact it is the intersection
of all Baer rings between $R$ and $Q(R)$. Moreover $B(R)$ can be constructed in the following way
$$B(R)=\{\sum_{i=1}^n r_ie_i \sep n\in\N, r\in R, e_i\in Q(R)\textrm{ with } e_i^2=e_i\}$$ 
\end{prop}

\begin{defi}
$B(R)$ as defined in the Proposition above is now usually known as the \emph{Baer hull of $R$}. Note that this should \emph{not be confused}
with $B(R)$ used by some authors (e.g. R.M. Raphael and A.S. Pierce) for the set of idempotents of $R$.
\end{defi}

\begin{lemma}
Let $A$ be a real closed ring and suppose $B\supset A$ for some ring $B$, then $A[e_{1},\dots,e_{n}]$ is a real closed ring for
$e_{i}\in B$ such that $e_{i}^{2}=e_{i}$ (i.e. $A$ adjoined with any idempotent of an over-ring is also a real closed ring).
\end{lemma}

\begin{proof}
We prove by induction. For $e\in B$ such that $e^{2}=e$, one can
immediately show that any element in $A[e]$ can be uniquely written
as $a_{1}e+a_{2}(1-e)$, thus we may write $A[e]\cong A_{1}\times A_{2}$
where $A_{1}$ is a commutative poring defined by
$$A_{1}:=\{ ae\sep a\in A\}$$
with multiplication operation $ae*be=abe$ and addition $ae+be=(a+b)e$,
similarly
$$A_{2}:=\{ a(1-e)\sep a\in A\}$$
with multiplication $a(1-e)*b(1-e)=ab(1-e)$ and addition $a(1-e)+b(1-e)=(a+b)(1-e)$.
Now the canonical maps $A\rightarrow A_{1}$ and $A\rightarrow A_{2}$
are surjective and the principal ideals $((e,0))$ and $((0,(1-e)))$
are radical in $A[e]$ thus we may regard $A_{1}$ and $A_{2}$ as
factor rings of $A$ modulo $((0,(1-e)))\cap A$ and $((e,0))\cap A$
(which are (real) radical ideals in our real closed rings $A$) respectively.
Since we know that the residue rings of a real closed rings modulo a real ideal
is real closed we also then know that $A_{1}$ and $A_{2}$ are real closed.
We also know that real closedness is stable upon taking arbitrary products, thus $A[e]$ is real 
closed. Now
if $A[e_{1},\dots,e_{n-1}]$ is real closed for idempotents $e_{1},\dots,e_{n-1}\in B$,
then take \[
A':=A[e_{1},\dots,e_{n-1}]\]
 and substitute it for $A$ in the argument above and use $e=e_{n}$.
\end{proof}

\begin{cor}\label{baercr}
Baer Hulls of real closed rings are real closed rings
as well.
\end{cor}

\begin{proof}
The Baer Hull $B(A)$ of a ring $A$ is constructed by adjoining it
with the idempotents of its complete ring of quotients, thus by taking
$$B(A):=\{\sum_{i=1}^{n}a_{i}e_{i}\sep n\in\N,a_{i}\in A,e_{i}^{2}=e_{i}\in Q(A)\}$$
 where $Q(A)$ is the complete ring of quotients of $A$. We know
from the previous Lemma that for any finite set of idempotents $e_{1},\dots,e_{n}\in Q(A)$
, $A[e_{1},\dots,e_{n}]$ is a real closed ring. But the family
$$\{ A[e_{1},\dots,e_{n}]\sep n\in\N,e_{i}^{2}=e_{i}\in Q(A)\}$$
 can be considered as a directed category with 
$$\sup(A[e_{1},\dots,e_{n}],A[f_{1},\dots,f_{m}])=A[e_{1},\dots,e_{n},f_{1},\dots,f_{m}]$$
and the canonical injection as the arrows. Having said that, one observes
that $B(A)$ is only the direct limit of objects of this category
(which is a subcategory of the real closed rings) in the category
of real closed rings. And since we know that the direct limits of real
closed rings are real closed we can conclude then that $B(A)$ is real closed.
\end{proof}

\begin{prop}\label{rcr_ic} Let $A$ be a real closed ring contained in a
real von Neumann regular ring $B$, then the integral closure
of $A$ in $B$ is\[
\bar{A}:=\{\sum_{i=1}^{n}a_{i}e_{i}\sep n\in\N,a_{i}\in A,e_{i}^{2}=e_{i}\in B\}\]
 and this is itself a real closed ring.
\end{prop}
\begin{proof}
Let $b\in B$ such that there exist a monic polynomial $f\in A[T]$
with $f(b)=0$. Denote $f_{\p}$ to be the projection of $f$ on $(A/(\p\cap A))[T]$,
for $\p\in\spec B$. Now $A/(\p\cap A)$ is also a real closed integral
domain, thus $A/(\p\cap A)$ is integrally closed in $B/\p$ (which
is a real closed field and contains $\qf{A/(\p\cap A)}$) and so we finally conclude that for any $\p\in\spec B$,
there exists an $a^{\p}\in A$ such that $a^{\p}(\p)=b(\p)$. 

Denote \[
K_{\p}:=\{\mathfrak{q}\in\spec B\sep(a^{\p}-b)(\mathfrak{q})=0\}=V(a^{\p}-b)\]
 we see that for each $\p\in\spec B$, $K_{\p}$ is nonempty and constructible
and the collection of these sets cover $\spec B$. Thus there are
$a_{1},\dots,a_{k}\in A$ such that for each $\p$ there is an $i\in\{1,\dots,k\}$
such that $a_{i}(\p)=b(\p)$. Also since $\spec B$ is von Neumann
regular (see \cite{Huckaba} Corollary 3.3(4)) there are idempotents $e_{i}$
such that $V(a_{i})=V(e_{i})$ in $\spec B$. Now define $e_{i}'$,
for $i=1,\dots,k$ in the following way : \begin{eqnarray*}
e_{1}' & = & e_{1}\\
e_{2}' & = & (1-e_{1}')e_{2}\\
e_{i}' & = & (1-\sum_{j=1}^{i-1}e_{j}')e_{i}\end{eqnarray*}
 then the claim is $b=\sum_{i=1}^{k}a_{i}e_{i}'$.

For each $i=1,\dots,k$ and $\p\in\spec B$ the reader can easily verify the following which also proves our claim
\begin{itemize}
\item $e_i'^2=e_i'$
\item If $e_i\in\p$ and $e_j\not\in\p$ for all $0\leq j< i$, then $e_i(\p)=1$ and 
$e_l(\p)=0$ for all $l=1,\dots,k$ with $l\neq j$
\end{itemize}
\end{proof}

\begin{theorem}\label{icofrcr} Let $A\subset B$ with $A,B$ being real closed rings. Then
the integral closure of $A$ in $B$ is a real closed ring. 
\end{theorem}
\begin{proof}
We can consider $A$ as a sub-poring of $\displaystyle\prod_{\sper B}\rho(\alpha)$.
Then by the Lemma above, we can construct the integral closure of
$A$ in $\displaystyle\prod_{\sper B}\rho(\alpha)$, this will be denoted as $\bar{A}$. We know that $\bar{A}$
is real closed. But the integral closure of $A$ in $B$ is in fact
$B\cap\bar{A}$ which can be considered as a fiber product of two
real closed rings, and we know that
fiber products of real closed rings are real closed as well.
\end{proof}

\begin{prop}\label{vnrcr}
A von Neumann regular ring is a real closed ring iff every residue field of it is a real closed field. In other words,
a von Neumann regular ring is a real closed ring iff for every $\p\in\spec A$, the field $A/\p$ is a real
closed field.
\end{prop}
\begin{proof}
This is just a result that can be automatically concluded from  Corollary 21 of \cite{rcr}
\end{proof}

\begin{defi}
Let $A$ be a ring and let $I\subset A$ be an ideal. We have the canonical residue map $\pi:A\rightarrow A/I$. Let now
$\bar f\in (A/I)[T]$ and say 
$$\bar f =\bar a_nT^n+\dots +\bar a_0$$
then there are some representatives $a_0,\dots,a_n\in A$ such that $\bar a_i=\pi(a_i)$. Then $f\in A[T]$ defined by 
$$f:=a_nT^n + \dots + a_0$$
is a \emph{lifting of $\bar f$ (in $A[T]$)} or that \emph{$\bar f$ is lifted to $f$} (obviously with this construction $f$ is not unique).
Conversely any polynomial $f\in A[T]$ is canonically mapped to $\bar f\in (A/I)[T]$ which we shall call the \emph{pushdown of $f$}
or we say \emph{$f$ is pushed down to $\bar f$ (in $(A/I)[T]$)}.
\end{defi}

\begin{lemma}
Let $K$ be a perfect field then every polynomial of odd degree in $K[T]$ has a zero in $K$ iff every polynomial with nonzero discriminant
of odd degree has a zero in $K$
\end{lemma}
\begin{proof}
One side of the equivalence is clear. The other side is also very easy, for any polynomial $f\in K[T]$ of odd degree can be written as 
$$f=\prod_{i=1}^n f_i$$
where $f_i$'s are irreducible for $i=1,\dots,n$ and without loss of generality we can assume $f_1$ has odd degree. 
Since our field is perfect $f_1$ is separable 
(i.e. has no multiple roots in the algebraic closure of $K$) and thus has a nonzero discriminant and therefore by assumption it has a 
zero in $K$ which is also a zero of $f$.
\end{proof}

%Jose
\begin{theorem}\label{vnricrcr}
Suppose $B$ is a real closed ring and $A$ a von Neumann regular sub-$f$-ring of $B$ which is integrally closed in $B$,
%is a von Neumann regular ring and $B$ is a real closed (von Neumann regular ring and an) over-poring of $A$. If $A$ is 
%integrally closed in $B$ 
then $A$ is also a real closed ring.
\end{theorem}
\begin{proof}
%By Lemma 1.14 of \cite{raphael} the canonical spectral map $\spec B\rightarrow \spec A$ is surjective.  
We make use of Proposition \ref{vnrcr} to prove our claim. We may let $\sum A^2$ to be the partial ordering of $A$ and consider 
$A$ a subporing of $B$. 
Choose any prime ideal $\p\in\spec A$, since $A$ is a von Neumann regular $f$-ring every prime ideal is an $l$-ideal 
(see for instance \cite{Capco} Proposition 3.7), thus by \cite{epi} Example 2.3 $A/\p$ is totally ordered by $\sum (A/\p)^2$.
If we show that for any (monic) polynomial $\bar f\in (A/\p)[T]$ of odd degree  $\bar f$ has a zero in $A/\p$ then we are done, since by \cite{KS}
Kapitel 1, Satz 1 $A/\p$ is a real closed field. 

Since all real fields are perfect fields, by the preceeding Lemma we need only deal with the case when $\bar f$ has a nonzero
discriminant. So let
$$\bar f(T):=T^{2n+1}+\sum_{k=0}^{2n} \bar a_kT^k \quad \bar a_k \in A/\p, k=1,\dots,2n$$ 
have a nonzero discriminant and suppose 
$$f(T):=T^{2n+1}+\sum_{k=0}^{2n} a_kT^k\quad  a_k \in A, k=1,\dots,2n$$
be a lifting of $\bar f$ in $A[T]$. 
%Note that the discriminant of $f\in A[T]$ is nonzero, since the discriminant of $\bar f$
%is the canonical image of the discriminant of $f$ in $A/\p$. 
Let $d\in A$ be the discriminant of $f$, then by \cite{PS} Theorem 2.1(v) $g\in A[T]$ defined by 
$$g(T):=T^{2n+1}+\sum_{k=0}^{2n} a_kd^{4n+2-2k}T^k$$
has a zero  $b\in B$, and since $A$ is integrally closed in $B$ we can conclude that $b\in A$. Let $\bar g\in (A/\p)[T]$ be the 
pushdown of $g$ in $(A/\p)[T]$, then the canonical image of $b$, say $\bar b$, in $A/\p$ is a zero of $\bar g$. 

Now the canonical image of $d$, say $\bar d$, in $A/\p$ is also the discriminant of $\bar f\in (A/\p)[T]$ and $\bar d\neq 0$. Since 
$A/\p$ is a field $\bar d$ has an inverse and so one easily checks that $\bar b\bar d^{-2}\in A/\p$ is a zero of $\bar f$.
\end{proof}

\begin{exam}
Now to show that there are real closed rings that are not real closed *. Consider the ring of continuous maps from $\R$ to $\R$
$$C(\R):=\{f:\R\rightarrow \R\sep f\textrm{ continuous }\}$$
this is  a real closed ring (see \cite{cont} Theorem 1.2) but is not Baer, because if 
$$f(x)=\left\{
\begin{array}{ll}
x  & x\geq 0\\
0	& x<0 \\
\end{array}
\right.$$
then $f$ has nontrivial annihilators which is strictly contained in our ring 
(for instance $f(-x)$ is an annihilator of $f$), but our ring has only trivial idempotents which clearly states
that our ring is not Baer.
\end{exam}

\comment{
As far as literature is concerned, it is not yet clear whether any ring $A$ has a unique real closure * (i.e. an essential 
integral extension of the ring which is real closed * and is isomorphic to any other extension of $A$ with such a property).
But we can at least show an interesting result if our ring $A$ is a real closed ring.
}

\begin{nota}
For any commutative reduced ring, by $Q(A)$ we mean the complete ring of quotients of $A$ and $B(A)$ as the Baer hull of 
$A$.
\end{nota}

Below is a Theorem which is probably known by many but I was not able to find a reference stating it. The technique of the proof
of the Theorem is a motivation from the proof of Proposition 2.2 in \cite{raphael} and in fact the proof is almost exactly the same.

\begin{theorem}\label{ieinee} Let $A$ be a reduced commutative ring, $B$ be an integral extension of $A$, and $C$ be an over-ring of $B$
and an essential extension of $A$. Then $B$ is an essential extension of $A$ (in words, any between integral extension
of a reduced ring and its essential extension is an essential extension).
\end{theorem}
\begin{proof} (Proof Technique from R.M. Raphael)
Let $b\in B\wo 0$ then $b\in C$ so $\exists c\in C$ such that $bc=:a\in A\wo 0$.
There also exists a monic polynomial $f\in A[T]$ with say 
$$f(T)=T^n + \sum_{i=0}^{n-1} a_{i}T^{i}\qquad n\in\N,a_i\in A\textrm{ for } i=1,\dots,n$$
such that $f(b)=0$.
% and $$n=\min\{\deg(g) \sep g\in A[T],g \textrm{ monic, and }g(b)=0\}$$
If $n=1$ there is nothing to prove so we are assuming here that $n\geq 2$.

Now if $a_0\neq 0$ then we are done since then 
$$b(b^{n-1}+\sum_{i=1}^{n-1} a_{i}b^{i-1})=-a_0\in A\wo{0}$$

Suppose therefore that $a_0=0$, we claim that there exists an $i\in\{1,\dots,n-1\}$ such that $a_ia^{i}\neq 0$. Suppose otherwise, then
$$c^nf(b)=c^n(b^n+\sum_{i=1}^{n-1} a_{i}b^{i})=a^n + \sum_{i=1}^{n-1} a_{i}a^{i}c^{n-i}=a^n=0$$
and because $A$ is reduced this implies that $a=0$ which is a contradiction.

Thus we can say that there exists an $m\in \N$ such that $m<n$ and such that $a^ma_m\neq 0$ and if $m\neq 0$ then 
$a^ma_i=0$ for $0\leq i<m$. More clearly stated 
$$m=\min\{j\sep a^ja_j \neq 0, j=1,\dots,n-1\}$$
So 
$$c^mf(b)=a^mb^{n-m}+a^m\sum_{i=m}^{n-1} a_{i}b^{i-m}=0$$
and thus $c^mf(b)-a^ma_m=-a^ma_m\in A\wo 0$ is a multiple of $b$ in $B$.
\end{proof}

An immediate Corollary to this Theorem follows

\begin{cor}\label{rcr-ic}
If $A$ is a real closed * ring and $B$ is an essential extension of $A$ then $A$ is integrally closed in $B$
\end{cor}

Here is a short but useful Lemma

\begin{lemma}\label{vnrideal}
If $A$ is a von Neumann regular ring, then every ideal is a radical ideal.
\end{lemma}
\begin{proof}
We need only show that every principal ideal is radical. If we show that for every $a\in A$ if $b^2\in aA$ then $b\in aA$, then the rest
of the proof will follow from induction. Suppose thus that $b^2\in aA$, and let $b'\in A$ be the pseudo-inverse of $b$, then 
$b'b^2=b\in aA$.
\end{proof}

It was hinted in \cite{rcvnr} Theorem 1 that von Neumann regular real closed * rings have real closed factor fields. 
In the proof of his Theorem, Zhizhong Dai made use of Theorem 10(i) in the work \cite{srcr}, in the latter Theorem 
however he stated and made use of the factor fields being real closed. Because of this confusion and because I could not see
a clear proof of this fact, I made an attempt and 
proved that the von Neumann regular real closed * rings have indeed real closed factor fields (using a result from \cite{raphael}).
But first we show a Lemma.

\begin{lemma}\label{baersubvnrf}
If $A$ is a Baer ring and $B$ a von Neumann regular $f$-ring which is an essential extension of $A$ then $A$ equipped with 
the partial ordering $A^+=B^+\cap A$ is an $f$-ring.
\end{lemma}
\begin{proof}
Given the condition for the partial ordering of $A$, one needs only check whether $a^+\in A$ for any $a\in A$.
According to \cite{raphael} Lemma 1.6, $A$ contains all the idempotents of $B$. Now it suffices to show that if 
$a\in A\subset B$ then $a^+\in B$ is actually in $A$. Let $(a^+)'$ be the pseudo-inverse of $a^+$ in $B$, then
$a^+(a^+)'\in A$ (because it is an idempotent in the von Neumann regular ring $B$). Thus
$$aa^+(a^+)'=(a^+-a^-)a^+(a^+)'=a^+a^+(a^+)'-a^-a^+(a^+)'=(a^+)^2(a^+)'-0=a^+\in A$$
\end{proof}

\begin{prop}\label{factorvnr}
If $A$ is a von Neumann regular real closed * ring then all its factor fields are real closed, in other words it is a von Neumann regular
real closed ring.
\end{prop}
\begin{proof}
The ring $A$ is Baer because of \cite{srcr} Theorem 7. So by the previous Lemma it suffices to show that $A$ has an essential extension
which is a von Neumann regular $f$-ring. The real closure $\rho(A)$ is also a von Neumann regular ring, this is clear by
the fact that the prime spectrum and the real spectrum are zero dimensional for a von Neumann regular ring and in \cite{SM} \S12 we see
that $\rho(A)$ has prime spectrum homeomorphic to the real spectrum of $A$ thus $\rho(A)$ is zero dimensional and reduced making it a 
von Neumann regular $f$-ring. If $\rho(A)$ were an essential extension of $A$ then we are done.

If $\rho(A)$ were not an essential extension of $A$ then there is a nonzero ideal $I$ of $\rho(A)$ such that $I\cap A={0}$, assume  $I$ is maximal
with the property that it is an ideal of $A$ and that $I\cap A={0}$. Then $\rho(A)/I$ is an essential extension of $A$, in fact $\rho(A)/I$
is also a von Neumann regular $f$-ring since any ideal of $\rho(A)$ is radical and thus any factor ring of $\rho(A)$ is a von Neumann regular
real closed ring (see Remark \ref{rcrprop}). Let $A^+$ be the original partial ordering of $A$. 
If we now equip $A$ with the partial ordering $A':=B^+\cap A$ we see by using Lemma \ref{baersubvnrf} that $(A,A')$ is an $f$-ring. 
(The claim is that $A'=A^+$.)

By Theorem \ref{ieinee} $A$ must
be integrally closed in $\rho(A)/I$ and thus by Theorem \ref{vnricrcr} $(A,A')$ is a real closed ring. But real closed rings can only have 
sum of squares as partial ordering (\cite{SM} \S12), and the partial orderings of reduced $f$-rings cannot be extended 
(See \cite{SM} Proposition 1.11). Thus $A^+=A'$ and $(A,A^+)$ is a real closed ring.
\end{proof}

\begin{nota}
If $R$ is a commutative ring, by $E(R)$ we mean the set of idempotents of $R$ i.e.
$$E(R):=\{e\in R\sep e^2 =e\}$$
with a ring structure (boolean ring) induced from the boolean algebra of $E(R)$, i.e. 
we define addition multiplication $*'$ and addition $+'$ of $e,f\in E(R)$ by
$$e+'f := e*(1-f)+f*(1-e) \quad e*'f=e*f$$
where $+,-,*$ are addition, subtraction and multiplication in $R$.
\end{nota}

The proof of Theorem 1 \cite{rcvnr} (apparently motivated by \cite{raphael}) is rather incomplete even if the result of 
the above Proposition is known. But if we compare with \cite{raphael} Proposition 2.7 (which has nothing to do with real rings)
one can derive a method of proof which make use makes use of a characterization in Satz 6.1 of \cite{Storrer} and the 
equivalence of reduced ringed spaces and the category of rings with conformal maps as developed by Pierce in \cite{Pierce}
Definition 2.1 and Theorem 6.6. This however hints us into delving in the sheaf structures of regular rings. Indeed upon consideration
of the sheaf structure of regular rings one derives an easier proof which can be applied to prove both Proposition 2.7 of \cite{raphael}
and Theorem 1 of \cite{rcvnr} (thus sparing us the effort of understanding the category of rings with conformal maps).

\comment{
We therefore attempt to write down the proof this Theorem (including the partial proof 
provided already by Z.Z. Dai). The method of proof will be very similar to \cite{raphael} Proposition 2.7.
}

%Jose
\begin{theorem} (\cite{rcvnr} Thereom 1)
A real von Neumann regular ring is real closed * iff it is Baer and all factor fields are real closed.
\end{theorem}
\begin{proof}
"$\Rightarrow$" The fact that a real closed * ring is Baer is rather clear, as it contains all the idempotents of its complete ring
 of quotients. That its factor fields are real closed come from the Proposition that we just proved.

"$\Leftarrow$" Assume that $R$ is a regular ring which is Baer and has real closed factor field and assume in addition that $S$ is a proper essential
 and integral extension of $R$ and is moreover a real ring.  Note first that by \cite{raphael} Corollary 1.10, $S$ will also turn out to be 
a von Neumann regular ring.

From Remark 1.17 of \cite{raphael} one realizes that $\spec S$ and $\spec R$ are canonically homeomorphic. Moreover $S$ is an integral 
extension of $R$, so for any $\p\in\spec S$ we have a canonical algebraic extension of real fields
$$R/(\p\cap R) \longrightarrow S/\p$$
this is an algebraic extension because $S$ was assumed to be integral over $R$. Now by assumption $R/(\p\cap R)$ is real closed thus 
it must be isomorphic to $S/\p$. Thus we find out that the factor fields of $S$ and that of $R$ are canonically isomorphic. Moreover since
we are dealing with von Neumann regular rings, the factor fields are equivalent to the stalk of the sheaf structure of the ring 
(i.e. if $A$ is a von Neumann regular ring then for any $\p\in\spec A$, $A/\p$ is isomorphic to $A_\p$, given by 
$a+\p\mapsto \frac{a}{1}$). Thus we gather that the sheaf structure of both $R$ and $S$ are isomorphic (the prime spectrum are isomorphic 
and the stalks are ring isomorphic making the sheaf isomorphic, see for instance \cite{Hart} Proposition 1.1). Thus so are the global
section of the sheaves (are ring isomorphic), which is none other than the rings $R$ and $S$.
\end{proof}

\begin{cor}\label{cqr_rcr=srcr}
If $A$ is a real closed ring, then $Q(A)$ is both real closed and real closed *
\end{cor}
\begin{proof}
That $Q(A)$ is real closed is a result in \cite{Gab} Corollary 6.11. Now because $Q(A)$ is a von Neumann regular ring 
which is obviously Baer and has real closed factor fields it is also real closed *.
\end{proof}

\begin{theorem}
A real closed ring $A$ is real closed * iff it is Baer
\end{theorem}
\begin{proof}
If any ring is real closed *, it is Baer (e.g. see \cite{srcr} Theorem 7). So suppose that we have a real 
closed ring
$A$ that is Baer, then we know at once from the previous Corollary that $Q(A)$ is real closed * and by 
the construction of the Baer Hull of $A$ in Proposition \ref{mborn} and the construction of the integral closure of $A$ in $Q(A)$ in
Proposition \ref{rcr_ic} we observe that $A$ is integrally closed in $Q(A)$. But a ring is real closed * iff it is integrally closed 
in its complete ring of quotient and that its complete ring of quotient is real closed * (see \cite{srcr} Theorem 3, it might have
been clear to Z.Z. Dai, but I did find the necessity to include Lemma 3.6 in \cite{raphael} to complete the last part of 
the proof of this Theorem by Dai).
\end{proof}

It is known that every real ring has a real closure * (i.e. a maximal integral extension that is real) and we know that this 
is not in general unique, however we can show uniqueness for the special case when our ring is real closed. In fact 
later on we make a more general statement that uniqueness can be expected if the extension $A\hookrightarrow\rho(A)$ is 
essential.

\begin{theorem}
A real closed ring $A$ has a unique real closure * which is ring isomorphic to $B(A)$ and thus also real closed.
\end{theorem}
\begin{proof}
Let $A$ be real closed and let $C$ be a real closure * of $A$. Note that because $B(A)$ is a real closed ring (see Corollary \ref{baercr}) and
because of the above Theorem, $B(A)$ is real closed *. Now $C$ is an essential integral extension of $A$, and by a result due to 
Storrer in Satz 10.1 \cite{Storrer} also stated in \cite{raphael} Theorem 3.12 we may regard $Q(A)$ as a subring of $Q(C)$. 
The injection $Q(A)\hookrightarrow Q(C)$ is also essential because the diagram below
$$
\begindc{\commdiag}[1]
\obj(0,60)[A]{$A$}
\obj(120,60)[B]{$Q(A)$}
\obj(0,0)[C]{$C$}
\obj(120,0)[D]{$Q(C)$}
\obj(60,60)[E]{$B(A)$}
\mor{A}{E}{}[1,3]
\mor{E}{B}{}[1,3]
\mor{B}{D}{}[1,3]
\mor{A}{C}{}[-1,3]
\mor{C}{D}{}[1,3]
\enddc
$$
is commutative and $A\hookrightarrow C\hookrightarrow Q(C)$
is essential.  We see that all the rings that we are dealing
with are subrings of $Q(C)$. $C$ is an integral extension of $A$ which is integrally closed in $Q(C)$, thus $C$ is the integral 
closure of $A$ in $Q(C)$ (because of the transitivity of integral extensions and because of Theorem 3 of \cite{srcr}).
We also observe, by Corollary \ref{rcr-ic}, that $B(A)$ must also be the integral closure of $A$ in $Q(C)$, since $B(A)$ is a real 
closed * and an integral extension of $A$ in $Q(C)$. Thus $B(A)$ and $C$ are isomorphic as rings.
\end{proof}

Thus we immediately have the following Corrolary

\begin{cor}
The real closure * of a real closed ring $A$ is a real closed ring.
\end{cor}

\comment{
\begin{proof}
Now let $C$ be a real closed * ring which is a maximal essential integral extension of $A$. By (ADD Reference)
$Q(C)$ is real closed * and $C$ is integrally closed in $Q(C)$. Since $Q(C)$ is a von Neumann regular Baer 
real closed * ring it is by Theorem 3 \cite{rcvnr} real closed. $C$ can also be considered as the integral 
closure of $A$ in $Q(C)$ (because $C$ is integrally closed in $Q(C)$ and it is an integral extension of $A$).
Thus by Theorem \ref{icofrcr} $C$ is real closed.
\end{proof}
}

We see thus that the real closure * of a real closed ring is a real closed ring. Now we try to determine 
the real closure 
of real closed * rings.

\begin{prop}
Let $A$ be a real closed * ring and suppose that $A\rightarrow\rho(A)$ is a flat or essential ring extension. Then $\rho(A)$ is Baer
and thus real closed *.
\end{prop}
\begin{proof}
Let $C$ be the smallest real closed ring between $A$ and $Q(A)$ (which is none other than the intersection of the real closed rings between 
$Q(A)$ and $A$).  By the choice of $C$, there is a unique 
surjective ring homomorphism $\rho(A)\twoheadrightarrow C$ that factors $A\hookrightarrow C$ (this is because the construction of
real closed ring results into a monoreflector from the category of porings to the category of real closed rings, see \cite{SM} \S12). 
We have the commutative diagram below 
$$
\begindc{\commdiag}[1]
\obj(0,60)[A]{$A$}
\obj(100,60)[B]{$C$}
\obj(0,0)[C]{$\rho(A)$}
\mor{A}{B}{}
\mor{A}{C}{}
%\mor{C}{B}{}
\arrow <1.5mm> [.2,1] from 10 10 to 90 50
\arrow <1.5mm> [.2,1] from 10 10 to 94 52
\enddc
$$
If $A\hookrightarrow \rho(A)$ were an essential extension then by the choice of $C$ and by the fact that $A\hookrightarrow C$ is essential,
the surjective ring homomorphism $\rho(A)\twoheadrightarrow C$ is actually a monomorphism and thus $C$ and $\rho(A)$ are ring isomorphic
having complete ring of quotients $Q(A)$ and being Baer. 

Now suppose that $A\hookrightarrow \rho(A)$ is a flat extension. Since $A$ is Baer it is known then that $Q(A)$ is a flat extension of $A$ 
(this can be found in say \cite{Findlay} Corollary 6.3).
Both $Q(A)$ and $\rho(A)$ being flat extensions of $A$ implies that the canonical morphisms
$\rho(A)\rightarrow \rho(A)\otimes Q(A)$ and $Q(A)\rightarrow \rho(A)\otimes Q(A)$ are monomorphisms of rings. 
%Set $B:=\rho(A)\otimes Q(A)$

We have the following commutative diagram
$$
\begindc{\commdiag}[1]
\obj(0,60)[A]{$A$}
\obj(120,60)[B]{$Q(A)$}
\obj(0,0)[C]{$\rho(A)$}
\obj(120,0)[D]{$Q(A)\otimes_A \rho(A)$}
\obj(60,60)[E]{$C$}
\mor{A}{E}{}[1,3]
\mor{E}{B}{}[1,3]
\mor{B}{D}{}[1,3]
\mor{A}{C}{}[-1,3]
\mor{C}{D}{}[1,3]
\enddc
$$
Note now that all our working rings are thus subrings of $Q(A)\otimes \rho(A)$. Because of the choice of $C$ and because 
real closed rings are closed under intersection, we must have 
$$C=C\cap\rho(A)$$
This clearly implies that $C\subset \rho(A)$ in $Q(A)\otimes \rho(A)$. We thus have the following
commutative diagram
$$
\begindc{\commdiag}[1]
\obj(0,60)[A]{$A$}
\obj(100,60)[B]{$C$}
\obj(0,0)[C]{$\rho(A)$}
\mor{A}{B}{}
\mor{A}{C}{}
%\mor{C}{B}{}
\arrow <1.5mm> [.2,1] from 10 10 to 90 50
\arrow <1.5mm> [.2,1] from 10 10 to 94 52
\arrow <1.5mm> [.2,1] from 94 46 to 14 6
\arrow <1.5mm> [.2,1] from 94 46 to 93.2 45.6
\enddc
$$
where $C\rightarrowtail \rho(A)$ is just the canonical injection. Now because $A\hookrightarrow \rho(A)$ is an epimorphic extension
(see \cite{SM} \S12, or \cite{Schw1} Chapter 2) in the category of 
reduced partially ordered rings (note that all our rings are ring homomorphisms are also poring morphisms, since $C$ is a real closed ring 
one can deduce
that for real closed ring the squares are the only partial ordering, the reader is referred to \cite{Schw1} Chapter 2 or \cite{SM} \S12)
, the monomorphism $A\hookrightarrow C$ is an essential epimorphic extension of $A$. And so because
the horizontal and vertical maps in the diagram above are epimorphism one
can then at once conclude that $\rho(A)\twoheadrightarrow C$ is an isomorphism. And finally we can say that $\rho(A)$ is Baer in both cases.
%$$0\rightarrow A\rightarrowtail\rho(A)\twoheadrightarrow C\rightarrow 0$$
%is an exact sequence in the category of rings. Moreover $
\end{proof}

\comment{
\begin{lemma}
If $Q(\rho(A))\cong\rho(Q(A))$ then $\rho(A)$ is Baer if $A$ is Baer.
\end{lemma}
}

\begin{theorem} Let $A$ be a reduced poring then
$\rho(A)$ is an essential extension of $A$ iff for all $a_1,\dots,a_n\in A$ such that 
$$\bar P(a_1,\dots,a_n)\subsetneq \sper A$$ 
there is an $a\in A$ such that 
$$\bar P(a_1,\dots,a_n)\subset Z(a)\subsetneq \sper A $$
\end{theorem}
\begin{proof}
Let $\tilde\rho:\sper \rho(A)\rightarrow \sper A$ be the spectral map of the injection 
$A\hookrightarrow \rho(A)$. We know from \cite{SM} §12 that $\tilde\rho$ and $\supp_{\rho(A)}$ are homeomorphisms.

"$\Rightarrow$" Let $a_1,\dots,a_n\in A$ and set 
$$X:=\bar P(a_1,\dots,a_n)$$
assuming that $X\subsetneq\sper A$. Then because $\tilde\rho$ is a homeomorphism and $X$ is a closed constructible set in $\sper A$, 
there is an $x\in\rho(A)$ such that $Z_{\rho(A)}(x)\subsetneq\sper\rho(A)$ and $\tilde\rho(Z_{\rho(A)}(x))=X$. Thus $x\neq 0$ 
(since for a reduced ring $B$, $x=0\Leftrightarrow Z(x)=\sper B$). Now since $\rho(A)$ is an essential extension of $A$, there must be
a $y\in\rho(A)$ such that $xy\in A\wo{0}$. Moreover one sees that 
$$Z_{\rho(A)}(x)\subset Z_{\rho(A)}(x)\cup Z_{\rho(A)}(y)=Z_{\rho(A)}(xy)$$ 
and if we set $a:=xy\in A$, we know from the fact that $\tilde\rho$ is a homeomorphism that 
$$\tilde\rho(Z_{\rho(A)}(a))=Z_A(a)\subsetneq\sper A$$
and therefore $X\subset Z_A(a)$.

"$\Leftarrow$" Let $x\in\rho(A)\wo{0}$ then $Z_{\rho(A)}(x)\subsetneq\sper\rho(A)$. Now because $Z_{\rho(A)}(x)$ is closed and constructible
in $\sper\rho(A)$ and $\tilde\rho$ is a homeomorphism
$$\exists\,a_1,\dots,a_n\in A \st \bar P_A(a_1,\dots,a_n)=\tilde\rho(Z_{\rho(A)}(x))\subsetneq\sper A$$
Set $X:=P_A(a_1,\dots,a_n)$ then we have an $a\in A\wo{0}$ such that $X\subset V_A(a)\subsetneq\sper A$. So
$$Z_{\rho(A)}(x)\subset\tilde\rho^{-1}(Z_A(a))=Z_{\rho(A)}(a)\subsetneq\sper\rho(A)$$
Now since $\supp_{\rho(A)}$ is a homeomorphism we have the identity 
$$\supp_{\rho(A)}(Z_{\rho(A)(x_0)}=V_{\rho(A)}(x_0)\quad\forall x_0\in\rho(A)$$
and thus
$$V_{\rho(A)}(x)\subset V_{\rho(A)}(a)\subset \spec\rho(A)$$
and common algebraic geometry knowledge gives us
$$\sqrt{\brac{x}}\supset\sqrt{\brac{a}}\Rightarrow \exists\,n\in\N,y\in\rho(A)\st a^n=xy\in A$$
Because $A$ is reduced we conclude that there exists $y\in\rho(A)$ such that $xy\in A\wo{0}$.
\end{proof}

%Jose
\begin{cor}
If $A$ is a real $f$-ring, then $\rho(A)$ is an essential extension of $A$
\end{cor}
\begin{proof}
We know from various sources that $f$-rings have $\supp$ which is a homeomorphism, making all the closed constructible 
sets of the form $Z(a)$  (see eg. \cite{Capco} Lemma 5.4 and \cite{epi} Example 2.3). Thus $f$-rings satisfy the condition of the 
Theorem above.
\end{proof}

%Jose
\begin{theorem}
If $A$ is any real ring such that $\rho(A)$ is an essential extension of $A$ then $A$ has a unique real closure * which is 
ring isomorphic to the integral closure of $A$ in $Q(\rho(A))$. In particular
real $f$-rings have unique real closure *.
\end{theorem}
\begin{proof}
Let $B$ be a real closure * of the ring $A$. Then by \cite{srcr} Theorem 3, $B$ is integrally closed in $Q(B)$ and $Q(B)$
is real closed * (and real closed, as these are equivalent for von Neumann regular rings). Thus we can consider the 
smallest real closed ring $C$ that lies between $A$ and $Q(B)$. Because of the choice of $C$ and because real closure 
construction forms a monoreflector we have the commutative diagram below 
$$
\begindc{\commdiag}[1]
\obj(0,60)[A]{$A$}
\obj(100,60)[B]{$C$}
\obj(0,0)[C]{$\rho(A)$}
\mor{A}{B}{}
\mor{A}{C}{}
%\mor{C}{B}{}
\arrow <1.5mm> [.2,1] from 10 10 to 90 50
\arrow <1.5mm> [.2,1] from 10 10 to 94 52
\enddc
$$
but then $\rho(A)\hookrightarrow C$ is a monomorphism because $\rho(A)$ is an essential extension of $A$, therefore $C$ is isomorphic to $\rho(A)$ . Thus $\rho(A)$ can always
be considered as a between ring of $A$ and the complete ring of quotient of any real closure * of $A$. Moreover for our ring $B$,
the extension $\rho(A)\hookrightarrow Q(B)$ is essential since the diagram below
$$
\begindc{\commdiag}[1]
\obj(0,60)[A]{$A$}
\obj(60,60)[B]{$B$}
\obj(0,0)[C]{$\rho(A)$}
\obj(120,60)[D]{$Q(B)$}
\mor{A}{B}{}[1,3]
\mor{A}{C}{}[1,3]
\mor{B}{D}{}[1,3]
\mor{C}{D}{}[1,3]
\enddc
$$
is commutative and the horizontal and vertical maps are essential extensions. Thus by \cite{Storrer} Satz 10.1, 
$Q(\rho(A))$ (which is real closed * by \cite{Gab} Corollary 6.11) can be considered as a subring of $Q(B)$. We observe
that all our working rings are subrings of $Q(B)$ and $B$ is the integral closure of $A$ in $Q(B)$ (by Theorem \ref{ieinee} 
and Corollary \ref{rcr-ic}). Analogously
as above we note that $Q(\rho(A))\hookrightarrow Q(B)$ is an essential extension. But since $Q(\rho(A))$ is real closed *, 
it is integrally closed in $Q(B)$ (see Corollary \ref{rcr-ic}). Thus the integral closure of $A$ in $Q(B)$ is none other
than the integral closure of $A$ in $Q(\rho(A))$.
\end{proof}

The following lemma is so easy that we will leave the proof to the reader.

\begin{lemma}
Let $I$ be some index set and let $A_i$ be integrally closed in $B_i$, with $A_i$ and $B_i$ being 
commutative unitary rings for all $i\in I$. Then $\prod_I A_i$ is integrally closed in $\prod_I B_i$.
\end{lemma}

%Jose
\begin{cor}
Direct products of real closed * rings are real closed.
\end{cor}
\begin{proof}
Let $\prod_I A_i$ be product of real closed * rings. Then by Utumi's Lemma (see \cite{lambek} 
Proposition 9 \S4.3) we get the following ring isomorphism
$$\prod_{i\in I} Q(A_i)\cong Q\prod_{i\in I} A_i$$
We also know that $A_i$ is integrally closed in $Q(A_i)$ for all $i\in I$, thus by the previous lemma
$\prod_I A_i$ is integrally closed in $\prod_I Q(A_i)$. Moreover by $\prod_I Q(A_i)$ is a von Neumann regular 
Baer real closed ring 
(products of real closed rings are real closed, see Remark \ref{rcrprop}) thus by \cite{rcvnr} it is a real closed * ring. Now by \cite{srcr} Theorem 3 a ring is 
real closed * iff it is integrally closed in its complete ring quotients and if its complete ring of quotients
is real closed *. Thus the proof is complete.
\end{proof}

Finally, below we produce an example of a real closed * ring whose real closure is neither an essential extension of the ring nor
is it a real closed * ring. Thus we end this paper by saying that the real closure of a real closed * ring is not necessarily 
real closed *.

%Josenow
\begin{exam}
Let $T$ be a real closed nonarchimedean field and $V$ be a subring of $T$ which is convex in $T$ (with respect to the partial ordering
of $T$) and thus a valuation ring of $T$. Suppose also in addition
that $V$ has at most 3 prime ideals, i.e. 
$$\spec V=\{\brac{0},\p,\mathfrak{q}\}$$
with $\brac{0}\subsetneq\p\subsetneq\mathfrak{q}$. Consider the ring $A:=S[x]+\p$ where $x\in\mathfrak q\backslash\p$ and 
$S$ is a representative field of $V$ i.e., a real closed field contained in $V$ and 
isomorphic to $V/\mathfrak{q}$ (see \cite{KS} Kapitel II \S5 Satz 3).

The ring $A$ is integrally closed in its quotient field $T$ which is a real closed ring and therefore by \cite{SV} 
Proposition 2 $A$ is real closed *. The prime spectrum and real spectrum of $\rho(A)$ are homeomorphic to $\sper A$, 
and one can compute the minimal prime cones
of the real spectrum of $A$. The minimal prime cones are of the form $\alpha+\alpha\p^2$, where 
$\alpha\in \sper S[x]$ are minimal. But the real and prime spectrums of polynomial rings over a real closed field are thoroughly 
studied and 
well known (see for instance examples found in \cite{KS} Kapitel III, p.110-111), it is also in general known that 
the space of the minimal prime cones of $S[x]$ do not form an extremally disconnected space (we can easily produce a basic 
open set consisting of Dedekind cuts whose closure is not open). Clearly the space of minimal prime cones of $A$ 
consist of more than one point, thus so does the minimal prime spectrum of $\rho(A)$ which tells us that $\rho(A)$ is 
not an integral domain. Had $\rho(A)$ been an essential extension over $A$, it would have been embedded between $A$ and 
$\qf A$ making it an integral domain, which is not the case. Since we concluded that the minimal prime cones of $\rho(A)$ 
(or $A$) are not extremally disconnected, this tells us that the space of minimal prime ideals of $\rho(A)$ is not extremally
disconnected which means that $\rho(A)$ cannot be Baer (see \cite{Mewborn} Proposition 2.3), thus nor can $\rho(A)$ be 
real closed *.
\comment{
(they are all 
Claim 1: If $\rho(A)$ were an essential extension of $A$ then it would be an integral domain \\
Claim 2: If $\rho(A)$ were a
}
\end{exam}

In a coming paper I will write about the uniqueness of real closure * for a wide set of rings and give examples of rings with 
different real closure *.

\begin{ack}
The author would like to thank Prof. Niels Schwartz and Marcus Tressl for their most valuable advises.
\end{ack}

\end{document}